\documentclass[10pt, titlepage]{amsart}

\usepackage{ae} 
\usepackage[T1]{fontenc}
\usepackage[cp1250]{inputenc}
\usepackage{amsmath}
\usepackage{amssymb, amsfonts,amscd,verbatim}
\usepackage[dvips]{graphicx}
\usepackage{latexsym}
\usepackage{indentfirst}
\usepackage{latexsym}

\usepackage{graphicx}
\usepackage{verbatim}   
\usepackage{color}      
\usepackage{subfigure}  

\usepackage{amsmath}    

\theoremstyle{plain}
\newtheorem{Prop}{Proposition}[section]
\newtheorem{Thm}[Prop]{Theorem}
\newtheorem{Cor}[Prop]{Corollary}

\theoremstyle{definition}
\newtheorem{Def}[Prop]{Definition}

\theoremstyle{remark}
\newtheorem{Rem}[Prop]{Remark}

\def\int{\mathop{\roman{int}}}

\def\1{^{-1}}

\def\R{{\mathbf R}}
\def\NN{{\mathbf N}}

\def\dim{\text{dim}}

\def\diam{\text{diam}}
\def\dist{\text{dist}}

\def\AAA{{\mathcal A}}

\def\EEE{{\mathcal E}}

\def\NNN{{\mathcal N}}

\def\UUU{{\mathcal U}}

\def\asdim{\mathrm{asdim}}
\def\An{\mathrm{An}}
\def\Lip{\mathrm{Lip}}
\def\OC{\mathrm{OpenCone}}
\def\dim{\mathrm{dim}}
\def\diam{\mathrm{diam}}
\def\dokaz{{\bf Proof. }}
\def\edokaz{\hfill $\blacksquare$}

\numberwithin{equation}{section}

\input pstricks.tex
\input xy
\xyoption{all}

\begin{document}
\title[
Sublinear Higson corona and Lipschitz extensions
]%
   {Sublinear Higson corona and Lipschitz extensions}

\author{M.~Cencelj}
\address{IMFM,
Univerza v Ljubljani,
Jadranska ulica 19,
SI-1111 Ljubljana,
Slovenija }
\email{matija.cencelj@guest.arnes.si}

\author{J.~Dydak}
\address{University of Tennessee, Knoxville, TN 37996, USA}
\email{dydak@math.utk.edu}

\author{J.~Smrekar}
\address{Fakulteta za Matematiko in Fiziko,
Univerza v Ljubljani,
Jadranska ulica 19,
SI-1111 Ljubljana,
Slovenija }
\email{jaka.smrekar@fmf.uni-lj.si}

\author{A.~Vavpeti\v c}
\address{Fakulteta za Matematiko in Fiziko,
Univerza v Ljubljani,
Jadranska ulica 19,
SI-1111 Ljubljana,
Slovenija }

\email{ales.vavpetic@fmf.uni-lj.si}

\date{ \today
}


\thanks{Supported in part by the Slovenian-USA research grant BI--US/05-06/002 and the ARRS
research project No. J1--6128--0101--04}
\thanks{The second-named author was partially supported
by Grant No.  2004047  from the United States-Israel Binational Science
Foundation (BSF),  Jerusalem, Israel.}

\begin{abstract}

The purpose of the paper is to characterize the dimension
of sublinear Higson corona $\nu_L(X)$ of $X$
in terms of Lipschitz extensions of functions:

\begin{Thm}
Suppose $(X,d)$ is a proper metric space.
The dimension of the sublinear Higson corona $\nu_L(X)$ of $X$
is the smallest integer $m\ge 0$ with the following property:
Any norm-preserving asymptotically Lipschitz function $f'\colon A\to \R^{m+1}$, $A\subset X$, extends to a 
norm-preserving asymptotically Lipschitz
function $g'\colon X\to \R^{m+1}$.
\end{Thm}

One should compare it to the result of Dranishnikov \cite{Dr1} who characterized the dimension of
the Higson corona $\nu(X)$ of $X$ is
 the smallest integer $n\ge 0$ such that $\R^{n+1}$
 is an absolute extensor of $X$ in the asymptotic category $\AAA$ (that means any proper asymptotically
Lipschitz function $f\colon A\to \R^{n+1}$, $A$ closed in $X$,
extends to a proper asymptotically Lipschitz function $f'\colon X\to \R^{n+1}$).
\par
 In \cite{Dr1} Dranishnikov introduced the category $\tilde \AAA$ whose objects
 are pointed proper metric spaces $X$ and morphisms are asymptotically Lipschitz
 functions $f\colon X\to Y$ such that there are constants $b,c > 0$ satisfying
 $|f(x)|\ge c\cdot |x|-b$ for all $x\in X$.
 We show $\dim(\nu_L(X))\leq n$
 if and only if $\R^{n+1}$ is an absolute extensor of $X$ in the category $\tilde\AAA$.
\par
As an application we reprove the following result
of Dranishnikov and Smith \cite{DRS}:

\begin{Thm}
Suppose $(X,d)$ is a proper metric space of finite asymptotic Assouad-Nagata
dimension $\asdim_{AN}(X)$. If $X$ is cocompact and connected, then $\asdim_{AN}(X)$
equals the dimension of the sublinear Higson corona $\nu_L(X)$ of $X$.
\end{Thm}
\end{abstract}

\maketitle

\medskip
\medskip
\tableofcontents

\section{Introduction}

In \cite{DRS} Dranishnikov and Smith related the dimension $\dim(\nu_L(X))$ of the sublinear Higson corona $\nu_L(X)$ of $X$ (see Section 4) to the asymptotic Assouad-Nagata dimension
$\asdim_{AN}(X)$ of $X$ as follows:

\begin{Thm}[Dranishnikov-Smith]\label{DRSThm}
Suppose $(X,d)$ is a connected cocompact proper metric space.
If $\asdim_{AN}(X)$ is finite, then $\asdim_{AN}(X)=\dim(\nu_L(X))$.
\end{Thm}

Since Assouad-Nagata dimension of $X$ is related to the sphere $S^m$
being a Lipschitz extensor of $X$ (see \cite{BDHM} and \cite{LangSch}), a natural problem arises of expressing $\dim(\nu_L(X))$
in terms of Lipschitz extensions.

Dranishnikov \cite{Dr1} (p.1105) introduced dimension $\dim^c(X)$ for proper
metric spaces $X$
as the smallest integer $n\ge 0$ such that any proper asymptotically
Lipschitz function $f\colon A\to \R^{n+1}$, $A$ closed in $X$,
extends to a proper asymptotically Lipschitz function $f'\colon X\to \R^{n+1}$.
Thus, $\dim^c(X)\leq n$ if and only if $\R^{n+1}$ is an absolute extensor
of $X$ in the asymptotic category $\AAA$ of \cite{Dr1}.
His basic result is that $\dim^c(X)$ equals the dimension of
Higson corona $\nu(X)$ of $X$ (\cite{Dr1}, Theorem 6.6 on p.1111).

\par In this paper we consider the case of a special class of
proper functions to $\R^{n+1}$: the norm-preserving functions.
Our spaces are assumed to be equipped with a base-point $x_0$.
In case of euclidean space $\R^m$ the base-point is always $0$.
Now, the {\it norm} $|x|$ of $x$ is defined to be $d(x,x_0)$.
$f\colon X\to \R^m$ is {\it norm-preserving} if $|f(x)|=|x|$ for all $x\in X$.

\par In analogy to $\dim^c(X)$ of Dranishnikov
we introduce the {\it sublinear asymptotic dimension} $\asdim_L(X)$
as the smallest integer $n\ge 0$ such that any norm-preserving asymptotically
Lipschitz function $f\colon A\subset X\to \R^{n+1}$
extends to a norm-preserving asymptotically Lipschitz function $f'\colon X\to \R^{n+1}$.
Our main result is that $\asdim_L(X)$ equals the dimension of
sublinear Higson corona $\nu_L(X)$ of $X$ if $X$ is proper.
\par Another way of generalizing $\dim^c$ would be to
look at the smallest integer $n\ge 0$ such that $\R^{n+1}$ is an absolute extensor of $X$
in the category $\tilde\AAA$ of \cite{Dr1}. In Section 6 we show that definition
to be equivalent to $\asdim_L$.
\par
As an application of our results we prove a slight generalization of \ref{DRSThm}
(see \ref{DimOfHigsonSublinearCoronaAndAN}).

\section{Norm-preserving Lipschitz functions}

An {\it $\epsilon$-net} in a metric space $X$ is a subset $X_1$ of $X$
such that $d(x,X_1) < \epsilon$ for all $x\in X$.
\par A metric space $X$ is {\it discrete} if there $\epsilon > 0$ such that
$d(x,y)\ge \epsilon$ for all $x,y\in X$, $x\ne y$. In that case the term
{\it $\epsilon$-discrete} will be used.
\par A function $f\colon (X,d_X)\to (Y,d_Y)$ of metric spaces
is {\it asymptotically Lipschitz} (or {\it large scale Lipschitz})
if there are numbers $\lambda, M\ge 0$ such that
$d_Y(f(x),f(y))\leq \lambda\cdot d_X(x,y)+M$ for all $x,y\in X$.
We will use the shortcut of $(\lambda,M)$-asymptotically Lipschitz in such case.
If $M=0$, then such maps are called {\it Lipschitz} or {\it $\lambda$-Lipschitz}.
Also, the notation $\Lip(f)\leq \lambda$ is being used in that case.
More precisely, $\Lip(f)$ is defined as $\sup\{\frac{d_Y(f(x),f(y))}{d_X(x,y)}\mid x\ne y\in X\}$.

\par Notice that if $X$ is $\epsilon$-discrete and $f\colon X\to Y$
is $(\lambda,M)$-asymptotically Lipschitz, then it is
$(\lambda+\frac{M}{\epsilon})$-Lipschitz.

\par Let us show how to relate $\asdim_L(X)$ to existence of Lipschitz
extensions.
\begin{Prop}\label{EquivOfLipAndAsympLipCase}
Suppose $(X,d)$ is a metric space and $X_1$ is a discrete $\epsilon$-net in $X$ for some $\epsilon > 0$.
If $m\ge 0$, then the following conditions are equivalent:
\begin{itemize}
\item[a.] Any asymptotically Lipschitz norm-preserving function $f'\colon A\subset X\to R^{m+1}$ 
 extends to an asymptotically Lipschitz norm-preserving function $F'\colon X\to R^{m+1}$.
\item[b.] Any Lipschitz norm-preserving function $f'\colon A\subset X_1\to R^{m+1}$ 
 extends to a Lipschitz norm-preserving function $F'\colon X_1\to R^{m+1}$.
\end{itemize}
\end{Prop}
\dokaz 
a)$\implies$b). Suppose $f'\colon A\subset X_1\to R^{m+1}$
is a norm-preserving Lipschitz function. Extend it to
an asymptotically Lipschitz function $F'\colon X_1\to R^{m+1}$
and notice that it is Lipschitz as $X_1$ is discrete.
\par b)$\implies$a).
Suppose $f'\colon A\subset X\to R^{m+1}$
is a norm-preserving $(\lambda,M)$-asymptotically Lipschitz function.
We may assume $x_0\in A$ and $f'$ is induced by $f\colon A\to S^m$
in the sense of $f'(x)=|x|\cdot f(x)$ for all $x\in A$.
Let $A_1$ be the set of all points $x$ in $X_1$ such that $d(x,A)<\epsilon$.
For each $x\in A_1\setminus A$ pick a point $a(x)\in A$
so that $d(x,a(x))< \epsilon$. If $x\in A_1\cap A$, put $a(x)=x$.
Define $g\colon A_1\to S^m$ by 
$g(x)=f(a(x))$.
We need the induced function $g'\colon A_1\to \R^{m+1}$ to be Lipschitz.
Given $x,y\in A_1$, we have $g'(x)-g'(y)=|x|\cdot f(a(x))-|y|\cdot f(a(y))=
f'(a(x))-f'(a(y))+(|x|-|a(x)|)\cdot f(a(x))+(|y|-|a(y)|)\cdot f(a(y))$,
so $|g'(x)-g'(y)|\leq \lambda d(a(x),a(y))+M+\epsilon+\epsilon\leq
\lambda\cdot (d(x,y)+2\epsilon)+M+2\epsilon$. Thus, $g'$ is asymptotically Lipschitz,
hence Lipschitz as $A_1$ is discrete.
Extend $g'$ over $X_1$ to a $\gamma$-Lipschitz norm-preserving
function $G'\colon X_1\to \R^{m+1}$. $G'$ is induced by
$G\colon X_1\to S^m$, an extension of $g$.
Given $x\in X\setminus (A\cup X_1)$ find a point $b(x)\in X_1$
such that $d(x,b(x)) < \epsilon$. For $x\in A\cup X_1$ put $b(x)=x$.
Extend $G$ to $F\colon X\to S^m$ by $F(x)=f(x)$ if $x\in A$ and $F(x)=G(b(x))$
if $x\notin A$.
Let $G'(x):=|x|\cdot G(x)$. Estimation of $|G'(x)-G'(y)|$ is obvious if both $x$ and $y$ belong to $A$. In case $x,y\notin A$ follow the argument above, so the only case of interest
is $x\in A$ and $y\notin A$.
Pick $x_1\in X_1$ satisfying $d(x,x_1)<\epsilon$ and notice $x_1\in A_1$, so $g(x_1)=f(a(x_1))$.
Now, $G'(x)-G'(y)=(G'(x)-G'(x_1))+(G'(b(y))-G'(y))+(G'(x_1)-G'(b(y))$
and $|G'(x_1)-G'(b(y))|\leq \gamma\cdot d(x_1,b(y))\leq \gamma\cdot (d(x,y)+2\epsilon)$,
$G'(b(y))-G'(y)=|b(y)|\cdot g(b(y))-|y|\cdot g(b(y))$ is of size at most $\epsilon$,
so the only remaining task is estimation of $G'(x)-G'(x_1)$.
However, $G'(x)-G'(x_1)=
f'(x)-f'(a(x_1))+(|x_1-|x|)\cdot f(x)+(|a(x_1)|-|x_1|)\cdot f(a(x_1))$ is of size at most $\lambda\cdot \epsilon+M+\epsilon+\epsilon$.
\hfill $\blacksquare$

In the remainder of the paper by $\An(X,r,s)$ we denote the annulus
$\{x\in X\mid r\leq |x| < s\}$, where $0\leq r \leq s$.

\begin{Prop}\label{BigClassOfSpecialHigsonFunctions}
Suppose $(X,d)$ is a discrete metric space, $r > 0$ and $M > 1$.
Put $X_k=\An(X,r\cdot M^{k-1},r\cdot M^{k+1})$ 
and $Y_k=\An(X, r\cdot M^{k-1},\infty)$.
If $f\colon X\to S^m$ is a function and $m\ge 0$, then the following conditions
are equivalent:
\begin{itemize}
\item[a.] $f$ induces
a norm-preserving $f'\colon X\to R^{m+1}$ that is Lipschitz.
\item[b.] The sequence $\{M^{k}\cdot \Lip(f\vert_{Y_k})\}_{k=1}^\infty$ is bounded.
\item[c.] The sequence $\{M^{k}\cdot \Lip(f\vert_{X_k})\}_{k=1}^\infty$ is bounded.
\end{itemize}
\end{Prop}
\dokaz a)$\implies$b). Suppose $f'$ is $\lambda$-Lipschitz.
Since $X$ is discrete, each $\Lip(f_k\vert_{Y_k})$ is finite.
If $x,y\in Y_k$ for some $k\ge 2$ and $|x|\leq |y|$, then
$\lambda\cdot d(x,y)\ge |f'(x)-f'(y)|=||x|f(x)-|y|f(y)|= | |x|(f(x)-f(y))+(|x|-|y|)f(y)|\ge
|x||f(x)-f(y)|-d(x,y)$,
so $M^k\frac{|f(x)-f(y)|}{d(x,y)}\leq M^k\frac{\lambda+1}{|x|}\leq M(\lambda+1)/r$.
\par c)$\implies$a). Suppose $\{M^{k}\cdot \Lip(f\vert_{X_k})\}_{k=1}^\infty$ is bounded by $C$.
Suppose $x,y\in X$ and $|x|\leq |y|$. If $x,y\in X_k$
for some $k$, then $||x|f(x)-|y|f(y)|= | |x|(f(x)-f(y))+(|x|-|y|)f(y)|\leq
|x||f(x)-f(y)|+d(x,y)\leq (|x|\cdot C/M^k+1)\cdot d(x,y)\leq (rMC+1)\cdot d(x,y)$.
If $x$ and $y$ do not belong to the same $X_k$,
then there is $n\ge 0$ such that $|x|\leq rM^n$ and $|y|\ge rM^{n+1}$.
In that case $d(x,y)\ge |y|-|x|\ge r(M-1)\cdot M^n\ge (M-1)|x|$
and $|f'(x)-f'(y)|=||x|f(x)-|y|f(y)|= | |x|(f(x)-f(y))+(|x|-|y|)f(y)|\leq
|x||f(x)-f(y)|+d(x,y)\leq 2|x|+d(x,y)\leq (2/(M-1)+1)\cdot d(x,y)$.
\hfill $\blacksquare$

\begin{Cor}\label{SuffCondForExtensionOfNormPresFunctions}
Suppose $(X,d)$ is a discrete metric space. If $m\ge 0$, then the following conditions are equivalent:
 \begin{itemize}
\item[a.] Every Lipschitz norm-preserving function $f'\colon A\subset X\to R^{m+1}$ 
 extends to a Lipschitz norm-preserving function $F'\colon X\to R^{m+1}$.
 \item[b.]  For any $r,s >0$, any $M > 1$,
 and any sequence
 of functions $f_k\colon A_k\subset \An(X,r\cdot M^{2k},r\cdot M^{2k+1})\to S^m$, $k\ge 1$, satisfying $\Lip(f_k)\leq \frac{s}{M^{2k}}$
 there is $c > 0$ and  functions $g_k\colon \An(X,r\cdot M^{2k},r\cdot M^{2k+1})\to S^m$ so that $\Lip(g_k)\leq \frac{c}{M^{2k}}$ and  $g_k\vert A_k=f_k$ for all $k\ge 1$.
  \item[c.]  For any $r,s >0$, any $M > 1$,
 and any sequence
 of functions $f_k\colon A_k\subset \An(X,r\cdot M^{k},r\cdot M^{k+1})\to S^m$, $k\ge 1$, satisfying $\Lip(f_k)\leq \frac{s}{M^{k}}$
 there is $c > 0$ and  functions $g_k\colon \An(X,r\cdot M^{k},r\cdot M^{k+1})\to S^m$ so that $\Lip(g_k)\leq \frac{c}{M^{k}}$ and  $g_k\vert A_k=f_k$ for all $k\ge 1$.
\end{itemize}
\end{Cor}
\dokaz (a)$\implies$(b). Let $A=\bigcup\limits_{k=1}^\infty A_k$.
Paste all $f_k$ to obtain $f\colon A\to S^m$. By
\ref{BigClassOfSpecialHigsonFunctions} the function $f'\colon A\to \R^{m+1}$
is Lipschitz, so it extends to a norm-preserving Lipschitz
$g'\colon X\to \R^{m+1}$ induced by $g$. Consider $g_k=g\vert_{\An(X,r\cdot M^{2k},r\cdot M^{2k+1})}$ and apply \ref{BigClassOfSpecialHigsonFunctions} again.
\par
(b)$\implies$(c).  It is a consequence of the fact that $\An(X,r\cdot M^{2k+1},r\cdot M^{2k+2})=\An(X,s\cdot M^{2k},s\cdot M^{2k+1})$, where $s=r\cdot M$,
 so one can apply b) twice and get c).
 \par (c)$\implies$(a).
 First observation is the following:
 \par 
For any $s >0$, any $M > 1$,
 and any sequence
 of functions $f_k\colon A_k\subset \An(X,M^{k},M^{k+3})\to S^m$ satisfying $\Lip(f_k)\leq \frac{s}{M^{k}}$
 there is $c > 0$ and  functions $g_k\colon \An(X,M^{k},M^{k+3})\to S^m$ so that $\Lip(g_k)\leq \frac{c}{M^{k}}$ and  $g_k\vert A_k=f_k$ for all $k\ge 1$.
\par It follows from setting $M'=M^3$ and applying c) three times
with values for $r$ being respectively $1$, $M$, and $M^2$.
\par Second observation is that
given $\lambda_i$-Lipschitz functions $u_i\colon X_i\to Y$, $i=1,2$,
from subsets of $X$ to a bounded metric space $Y$
such that $\dist(X_1\setminus X_2,X_2\setminus X_1)\ge \mu$
and $u_1\vert_{X_1\cap X_2)}=u_2\vert_{X_1\cap X_2)}$,
then $u=u_1\cup u_2\colon X_1\cup X_2\to Y$
is $\max(\lambda_1,\lambda_2,\frac{\diam(Y)}{\mu})$-Lipschitz.
\par
Suppose $f'\colon A\to R^{m+1}$ is a norm-preserving
Lipschitz function and $A\subset X$. Define $f(a)=f'(a)/|a|$
for $a\in A\setminus \{x_0\}$ and define $f(x_0)$ arbitrarily.
By \ref{BigClassOfSpecialHigsonFunctions} there is
$s > 0$ such that $\Lip(f\vert_{\An(A,2^k,\infty)})\leq s/2^k$ for all $k\ge 0$.
Find $t > 0$ and $f_k\colon \An(X,2^{k-1},2^{k+2})\to S^m$
for each $k\ge 2$ so that $f_k$ extends $f\vert_{\An(A,2^{k-1},2^{k+2})}$
and $\{2^k\cdot \Lip(f_k)\}$ is bounded.
Paste $f_k$ and $f_{k+2}$ to get $g_k$ on $\An(X,2^k,2^{k+1})\cup \An(A,2^{k+1},2^{k+2})\cup \An(X,2^{k+2},2^{k+3})$ so that $\{2^k\cdot \Lip(g_k)\}$ is bounded.
Get extensions $h_k\colon  \An(X,2^{k},2^{k+3})\to S^m$ of $g_k$
so that $\{2^k\cdot \Lip(h_k)\}$ is bounded. Splice even-numbered $h_k$ to get
$h\colon X\to S^m$ so that $h'$ is Lipschitz and $h$ extends $f$.
\hfill $\blacksquare$

\begin{Cor}\label{UpperHemisphereCase}
Let $(X,d)$ be a discrete metric space and $m\ge 0$.
If $K\subset S^m$ is a Lipschitz extensor of $X$
and $f\colon A\subset X\to K$ induces Lipschitz function
$f'\colon X\to \R^{m+1}$, then it extends to
$g\colon X\to K$ inducing Lipschitz
function $g'\colon X\to \R^{m+1}$.
\end{Cor}
\dokaz 
$K$ being a Lipschitz extensor of $X$
means that there is $C > 0$ such that any
$\lambda$-Lipschitz $f\colon A\subset X\to K$
extends over $X$ to a $C\cdot \lambda$-Lipschitz function.
Follow the same procedure as in the proof of
\ref{SuffCondForExtensionOfNormPresFunctions} using the fact
$K$ is a Lipschitz extensor of $X$.
\hfill $\blacksquare$

\begin{Cor}\label{IfSmIsLipschitzExtensorThenDimIsAtMostm}
Suppose $(X,d)$ is a discrete metric space.
If $S^m$ is a Lipschitz extensor of $X$, then
any norm-preserving Lipschitz function $f\colon A\to \R^{m+1}$ extends to a 
norm-preserving Lipschitz
function $f'\colon X\to \R^{m+1}$.
\end{Cor}

In \cite{DRS} the concept of a {\it cocompact} metric space $X$
was introduced to mean that there is a compact subset $K$
of $X$ such that for any $x\in X$ there is an isometry $\gamma$
of $X$ with $x\in \gamma(K)$. In other words,
the group of isometries of $X$ acts on $X$ cocompactly.
Analogously, one can introduce the concept of a {\it cobounded} space
$X$ to mean that Isom$(X)$ acts on $X$ coboundedly,
i.e. there is a bounded subset $B$ of $X$
with $X=\bigcup\limits_{\gamma\in Isom(X)}\gamma(B)$.

\par $(X,d)$ is {\it $M$-scale connected} if any pair of points $x,y\in X$
can be connected by a chain $x_1=x,\ldots, x_k=y$ such that $d(x_i,x_{i+1})\leq M$
for all $i\leq k-1$.

\begin{Cor}\label{SexySuffCondForExtensionOfNormPresFunctions}
Suppose $(X,d)$ is a metric space and $\epsilon > 0$.
If $X$ is cobounded and $M$-scale connected for some $M > 0$, then the following conditions are equivalent:
\begin{itemize}
\item[a.] Any asymptotically Lipschitz norm-preserving function $f'\colon A\subset X\to R^{m+1}$ 
 extends to an asymptotically Lipschitz norm-preserving function $F'\colon X\to R^{m+1}$.
\item[b.] There is a function $c\colon \R_+\to\R_+$
 such that any $f\colon B\subset Y\to S^m$, $Y$ an $\epsilon$-discrete subset of $X$, satisfying $\Lip(f)\leq \frac{s}{\diam(Y)}$ for some $s> 0$
 extends to $g\colon Y\to S^m$ so that $\Lip(g)\leq c(s)\cdot \Lip(f)$.
\end{itemize}
\end{Cor}
\dokaz a)$\implies$b). Suppose there is $s > 0$ such that for each $n\ge 1$
one can find $f_n\colon B_n\subset Y_n\to S^m$, $Y_n$ an $\epsilon$-discrete subset of $X$, satisfying $\Lip(f_n)\leq \frac{s}{\diam(Y_n)}$
but any extension $g\colon Y_n\to S^m$ of $f_n$ satisfies $\Lip(g) > \alpha(n)\cdot \Lip(f_n)$,
where $\alpha(n)\to\infty$ as $n\to\infty$.
Notice that each $f_n$ has an extension $g$ such that $\Lip(g)\leq \frac{2}{\epsilon}$.
That means $\Lip(f_n)\leq \frac{2}{\epsilon\cdot \alpha(n)}$ and $\Lip(f_n)\to 0$. 
By enlarging $M$ if necessary we may assume that $M > 6$ and
each $\An(X,2M^k,\frac{1}{3}M^{k+1})$ contains a point $x_k$. For each $n\ge 1$ pick
$k(n)$ so that $M^{k(n)}\leq \frac{2Ms}{\Lip(f_n)} < M^{k(n)+1}$
and by applying an isometry assume $d(x_{k(n)},Y_n) < M$ for all $n\ge 1$.
We may assume $k(n)$ is a strictly increasing function and $k(n)\ge 2$ for all $n\ge 1$ (by switching to
a subsequence if necessary).
Notice $\diam(Y_n)\leq \frac{s}{\Lip(f_n)}\leq M^{k(n)}/2$.
Given $x\in Y_n$ one has $|x_{k(n)}|-M-\diam(Y_n) < |x| < |x_{k(n)}|+M+\diam(Y_n)$,
so $2M^{k(n)}-M-M^{k(n)}/2\leq |x| < \frac{1}{3}M^{k(n)+1}+M+M^{k(n)}/2$ and $M^{k(n)}\leq |x| < M^{k(n)+1}$ as $M > 6$.
That means $Y_n\subset \An(X,M^{k(n)},M^{k(n)+1})$.

Since $\Lip(f_n)\leq \frac{2Ms}{M^{k(n)}}$ for all $n\ge 1$, \ref{SuffCondForExtensionOfNormPresFunctions} implies existence of $c > 0$ and extensions $g_n\colon Y_n\to S^m$ such that $\Lip(g_n)\leq \frac{c}{M^{k(n)}}$.
Thus, for all $n\ge 1$, $\frac{c}{M^{k(n)}}\ge \Lip(g_n)\ge \alpha(n)\cdot \Lip(f_n)\ge \alpha(n)\cdot \frac{2Ms}{M^{k(n)}}$ resulting in $\frac{c}{2Ms}\ge \alpha(n)$ for all $n$,
a contradiction.
\par b)$\implies$a). Use \ref{SuffCondForExtensionOfNormPresFunctions}
and \ref{EquivOfLipAndAsympLipCase}.
\hfill $\blacksquare$

\begin{Cor}\label{CharOfAsympANDimViaNormPres}
Suppose $(X,d)$ is a metric space such that
any asymptotically Lipschitz norm-preserving function $f'\colon A\subset X\to R^{m+1}$ 
 extends to an asymptotically Lipschitz norm-preserving function $F'\colon X\to R^{m+1}$
 for some $m\ge 0$.
If $X$ is cobounded, $M$-scale connected for some $M > 0$, and $\asdim_{AN}(X)\leq m+1$, then 
$\asdim_{AN}(X)\leq m$.
\end{Cor}
\dokaz
Suppose any asymptotically Lipschitz norm-preserving function $f'\colon A\subset X\to R^{m+1}$ 
 extends to an asymptotically Lipschitz norm-preserving function $F'\colon X\to R^{m+1}$.
 Pick a $1$-net $X_1$ in $X$.
By \ref{SexySuffCondForExtensionOfNormPresFunctions} there is a function $c\colon \R_+\to\R_+$
 such that any $f\colon B\subset Y\to S^m$, $Y\subset X_1$, satisfying $\Lip(f)\leq \frac{s}{\diam(Y)}$
 extends to $g\colon Y\to S^m$ so that $\Lip(g)\leq c(s)\cdot \Lip(f)$.
 As Assouad-Nagata dimension of $X_1$ is at most $m+1$,
 there is a constant $C > 0$ such that for each $r > 0$
 there is a cover $\UUU^r$ of $X_1$ that can be expressed as the union
 of $m+2$ families $\UUU^r_i$, $1\leq i\leq m+2$, so that each $\UUU_i^r$
 is $r$-disjoint and the diameter of elements of $\UUU^r$ is at most $C\cdot r$.
 Let $A_r$ be all the points of $X_1$ contained in at most $m+1$ elements of
 $\UUU^r$.
 Consider the canonical map $\phi^r\colon X_1\to \NNN(\UUU^r)$
 to the nerve $\NNN(\UUU^r)$ of $\UUU^r$ for each $r > 0$.
 Notice $\Lip(\phi^r)\leq \frac{\lambda}{r}$ for some constant
 $\lambda > 0$
 (see
\cite{BD} and \cite{BS}) and $\diam((\phi^r)^{-1})(\Delta))\leq 2Cr$
 for each $r > 0$ and each $(m+1)$-simplex $\Delta$ of $\NNN(\UUU^r)$.
 Let $X^r_\Delta=(\phi^r)^{-1}(\Delta)$ and $\phi^r_\Delta=\phi^r\vert_{X^r_\Delta}\colon X^r_\Delta\to\Delta$ for any $(m+1)$-simplex $\Delta$ of $\NNN(\UUU^r)$.
 Thus $\Lip(\phi^r_\Delta)\leq \frac{2\lambda}{\diam(X^r_\Delta)}$ 
 and for each $r > 0$ and for each $(m+1)$-simplex $\Delta$ of $\NNN(\UUU^r)$
 there is an extension $\psi^r_\Delta\colon X^r_\Delta\to \partial\Delta$
 of $\phi^r_\Delta\vert (X^r_\Delta\cap A_r)$ satisfying $\Lip(\psi^r_\Delta)\leq t\cdot \Lip(\phi^r_\Delta)$, where $t=c(2\lambda)$.
 Given a vertex $U$ of $\UUU^r$ consider all $\Delta$ containing $U$ as a vertex
 and let $s(U)$ be the union of all $(\psi^r_\Delta)^{-1}(st(U))$ and $U\cap A_r$. That means $s(U)$ shrinks $U$ to points that are either
 contained in $A_r\cap U$ or their $U$-barycentric cooridinate
 of $\psi^r$ is positive. Thus the mesh of the cover $\{s(U)\}_{U\in\UUU^r}$ is at most $C\cdot r$, so it suffices to show that its Lebesque number
 is at least $K\cdot r$ for some constant $K$ (independent on $r$)
 and the multiplicity of that cover is at most $m+1$.
 \par The reason multiplicity does not exceed $m+1$ is that
 each $x\in X$ is being mapped to the $n$-skeleton of the nerve
 by maps $\psi^r_\Delta$.
 \par
 Given $x\in X$ find all elements $U_1,\ldots,U_k$ of $\UUU^r$
 containing it. If $k< m+2$, then the ball $B(x,r/2)$
 intersects exactly one $U_i$ and $s(U_i)$ must contain $B(x,r/2)\subset A_r\cap U_i$. If $k=m+2$, then for $\Delta=[U_1,\ldots,U_{m+2}]$
 pick $i$ such that $U_i$-barycentric coordinate of $\psi^r_\Delta(x)$
 is at least $\frac{1}{m+2}$.
 If $y\in B(x,\frac{r}{\lambda(m+2)t})$ and the $U_i$-barycentric coordinate
 of of $\psi^r_\Delta(y)$ is $0$,
 then $1/(m+2)\leq |\psi^r_\Delta(x)-\psi^r_\Delta(y)|\leq d(x,y)\cdot t \lambda/r<  1/(m+2)$, a contradiction.
 Thus $B(x,\frac{r}{\lambda(m+2)t})\subset s(U_i)$
 and the Lebesque number of $\{s(U)\}_{U\in\UUU^r}$ is at least $K\cdot r$,
 where $K=\frac{1}{\lambda(m+2)t}$.
 That means $m\ge \dim_{AN}(X_1)=\asdim_{AN}(X_1)=\asdim_{AN}(X)$.
\hfill $\blacksquare$

\section{Open cones and partitions of unity}

Consider a compact subset $K$ of $\R^m\setminus \{0\}$ such that
each ray from $0$ intersects at most one point of $K$.
The union of rays from $0$ via points of $K$ is denoted by $\OC(K)$.
Thus, every point of $\OC(K)\setminus \{0\}$ has unique representation
as $t\cdot k$, $t > 0$ and $k\in K$. The most important examples
are $K=S^{m-1}$ in which case $\OC(K)=R^m$
and $K\subset \Delta^{m-1}$, $\Delta^{m-1}$ being the standard
simplex spanned by the basis $e_i$, $1\leq i\leq m$, of $\R^m$.

We are interested in open cones over sets $K$ such that any point
$x\in K$ has a closed neighborhood $A_x$ in $K$ with the property that
the distance between any \lq tangent\rq\ vector $(v-w)/|v-w|$, $v,w\in A_x$,
and $x/|x|$ is bigger than $\alpha(x)$,
where $\alpha\colon K\to (0,\infty)$ is a function.
For simplicity, lets call such sets {\it regular}. Notice that both
$S^{m-1}$
and $\Delta^{m-1}$ are regular.

\par The reason for next proposition is that simplexes
offer better geometrical properties (convexity) than subsets
of spheres. Our main use for it is in replacing $\R^{m+1}$
by the open cone over $\Delta^m$.

\begin{Prop}\label{ConeEquivalence}
Suppose $f\colon K\to L$ is a Lipschitz map of compact subsets
of euclidean spaces.
If $K$ is regular, then
$$f'\colon \OC(K)\to \OC(L)$$
 defined by
$f'(t\cdot k)=t\cdot f(k)$ is Lipschitz.
In particular, if $f\colon K\to L$ is a bi-Lipschitz homeomorphism of regular sets, then
$f'$ is a bi-Lipschitz homeomorphism.
\end{Prop}
\dokaz Suppose there is a sequence of points $t_n\cdot k_n, s_n\cdot l_n\in \OC(K)$
such that $|t_nf(k_n)-s_nf(l_n)| > \beta(n)|t_nk_n-s_nl_n|$
and $\beta(n)\to\infty$ as $n\to\infty$.
We may assume $s_n\ne 0$ (only one of $t_n,s_n$ can be $0$ for a given $n$).
Replacing $t_n$ by $t_n/s_n$ and $s_n$ by $1$ (divide both sides of inequality by $s_n$), we may assume $s_n=1$ for all
$n$. Notice $0$ cannot be the limit of an infinite subsequence of $\{t_n\}$
(in that case $\beta(n) < 2\diam(L)/d(0,K)$ for infinitely many $n$,
a contradiction). Similarly, no subsequence of $\{t_n\}$ can approach infinity
(reduce to the previous case by switching to $1/t_n$ via dividing).
Without loss of generality, we may assume $t_n\to t_0 > 0$
and $k_n\to k_0$, $l_n\to l_0$ (as $K$, $L$ are compact).
Notice $t_n=1$ only for finitely many $n$ (as $f$ is Lipschitz), so assume $t_n\ne 1$ for all $n$.
Similarly, we may assume $k_n\ne l_n$ for all $n$.
Now $t_0=1$ (otherwise the minimum distance from $K$ to $t_0K$ is positive
and one can get a bound on $\beta(n)$) which implies $k_0=l_0$.
It is here where we use regularity of $K$: assume all $k_n,l_n$ lie in a closed
neighborhood $A$ of $k_0=l_0$ so that there is $\alpha > 0$
with the property that $|(k-l)/|k-l|-tp| > \alpha$ for all $k,l,p\in A$ and all $t\in\R$
(the original assumption was satisfied for all $t > 0$ but it implies such inequality for all $t$).
Put $(t_n-1)/|k_n-l_n|=r_n$ and $(k_n-l_n)/|k_n-l_n|=v_n$.
The inequality $|t_nf(k_n)-s_nf(l_n)| > \beta(n)|t_nk_n-s_nl_n|$
can be easily transformed into $|r_nf(k_n)+(f(k_n)-f(l_n))/|k_n-l_n|| > \beta(n)|r_nk_n+v_n|$.
As $|f(k_n)-f(l_n)|/|k_n-l_n|$ is bounded ($f$ is Lipschitz),
no subsequence of $r_n$ can be unbounded, so we may assume $r_n\to r_0$.
That means $r_nk_n+v_n\to 0$ and $r_0=-1/|k_0|$.
In that case $|r_nk_n+v_n|> \alpha/2$ for large $n$, hence $|r_nf(k_n)+(f(k_n)-f(l_n))/|k_n-l_n||$ is bounded resulting in $\beta(n)$ being bounded, a contradiction.
\hfill $\blacksquare$

The next result shows benefits in switching to open cones over
simplexes. Namely, it provides a simple criterion
for the canonical partition of unity subordinate to a cover to induce a Lipschitz
\lq norm-preserving\rq\ map.

Given an open cover $\UUU=\{U_1,\ldots,U_k\}$ of a metric space $X$ by proper subsets,
by the {\it canonical partition of unity $\phi_\UUU$ induced by $\UUU$}
we mean functions $\phi_i(x):=\frac{d(x,X\setminus U_i)}{\sum\limits_{j=1}^k d(x,X\setminus U_j)}$, $i=1,\ldots,k$. Such a partition of unity
gives the {\it canonical map} from $X$ to the nerve $\NNN(\UUU)$
of $\UUU$ defined by $\phi(x)=\sum\limits_{i=1}^k \phi_i(x)\cdot v_i$,
where $v_i$ is the vertex corresponding to $U_i$.

\begin{Prop}\label{CharOfNiceCoversviaLip}
Suppose $\UUU=\{U_1,\ldots,U_k\}$ is an open cover of a metric space $X$
by proper subsets
and $\phi_\UUU$ is the canonical partition of unity induced by $\UUU$.
$\phi_\UUU'\colon X\to \OC(\NNN(\UUU))$
is Lipschitz if and only if 
there is $\epsilon > 0$ such that 
one has $$\sum\limits_{i=1}^k d(x,X\setminus U_i)\ge \epsilon |x|$$
for all $x$.
\end{Prop}
\dokaz Fix $M > 1$.
Choose $i\leq k$ and put $f(x)=d(x,X\setminus U_i)$,
$g(x)=\sum\limits_{i=j}^k d(x,X\setminus U_j)-f(x)$.
Notice $f$ and $g$ are $k$-Lipschitz and $f(x)+g(x)\ge \epsilon |x|$.
Let $h=\frac{f}{f+g}$. We want to show $|h(x)-h(y)|\leq \frac{3kd(x,y)}{\epsilon |y|}$
if $|x|\ge |y| > 0$ for some $x,y\in X$.
Put $a=\frac{3kd(x,y)}{\epsilon |y|}$. If $h(x)-h(y) > a$, then
$\frac{f(x)}{f(x)+g(x)}-\frac{f(x)-kd(x,y)}{f(x)+g(x)+2\cdot kd(x,y)} > a$ as well.
Since $\frac{f(x)}{f(x)+g(x)}-\frac{f(x)-kd(x,y)}{f(x)+g(x)+2\cdot kd(x,y)} =
\frac{f(x)\cdot 2kd(x,y)+kd(x,y)\cdot (f(x)+g(x))}{(f(x)+g(x))\cdot (f(x)+g(x)+2kd(x,y))}
\leq \frac{3kd(x,y)}{f(x)+g(x)+2kd(x,y)}\leq \frac{3kd(x,y)}{f(x)+g(x)}\leq a$,
we arrive at a contradiction.
\par Now, if $y=x_0$, then $||x|h(x)-|y|h(y)|=|x|h(x)\leq |x|=d(x,y)$.
If $|x|\ge |y| > 0$, then $||x|h(x)-|y|h(y)|=||y|(h(x)-h(y)+h(x)(|x|-|y|)|\leq
3kd(x,y)/\epsilon+d(x,y)=(3k/\epsilon+1)d(x,y)$.
\par Suppose $\phi'$ is $\lambda$-Lipschitz.
Let $S(x)=\sum\limits_{i=1}^k d(x,X\setminus U_i)$.
Suppose $x\in X$ and choose $i$ such that $x\in U_i$.
Find $y\in X\setminus U_i$ so that $d(x,X\setminus U_i)> d(x,y)/2$.
Since $\lambda\cdot d(x,y)S(x)\ge |x|d(x,X\setminus U_i)-|y|d(y,X\setminus U_i)S(x)/S(y)=|x|d(x,X\setminus U_i)> |x|d(x,y)/2$, we get $S(x) > \frac{|x|}{2\lambda}$.
\hfill $\blacksquare$

Our last result on open cones shows that functions
from $X$ to $\Delta^m$ that induce Lipschitz maps from $X$ to $\OC(\Delta^m)$
form a \lq convex\rq\ set.

\begin{Prop}\label{ConvexityOfSpecialHigson}
Suppose $(X,d)$ is a metric space and $\gamma\colon X\to \Delta^1$ 
induces a Lipschitz function $\gamma'\colon X\to \OC(\Delta^1)$.
If $f,g\colon X\to \Delta^m$ induce Lipschitz
functions $f',g'\colon X\to \OC(\Delta^m)$ for some $m\ge 1$,
then $h=\alpha\cdot f+\beta\cdot g\colon X\to \Delta^m$,
where $\alpha,\beta$ are barycentric coordinates of $\gamma$,
induces a Lipschitz
function $h'\colon X\to \OC(\Delta^m)$.
\end{Prop}
\dokaz 
$h'(x)-h'(y)=|x|\alpha(x)f(x)+|x|\beta(x)g(x)-|y|\alpha(y)f(y)-|y|\beta(y)g(y)=
\alpha(x)(f'(x)-f'(y))+(\alpha(x)-\alpha(y))|y|f(y)+
\beta(x)(g'(x)-g'(y))+(\beta(x)-\beta(y))|y|g(y)$.
Also, $|(\alpha(x)-\alpha(y))\cdot |y||=|(|x|\alpha(x)-|y|\alpha(y))+\alpha(x)\cdot (|y|-|x|)|\leq
|\gamma'(x)-\gamma'(y)|+d(x,y)\leq (1+\Lip(\gamma'))\cdot d(x,y)$
(the same holds for $\beta$),
so $|h'(x)-h'(y)|\leq \Lip(f')\cdot d(x,y)+\Lip(g')\cdot d(x,y)+(2\Lip(\gamma')+2)\cdot d(x,y)$
and $h'$ is Lipschitz.
\hfill $\blacksquare$

\section{Sublinear coarse structure}

The sublinear coarse structure was introduced in \cite{DRS}. Our approach is a little bit different to allow for a better transition from previous sections of our paper.

\begin{Def}\label{SublinearMapsDef}
A continuous function $s\colon \R_+\to\R_+$ is called {\it asymptotically
sublinear} if for each non-constant linear function $f\colon \R_+\to\R_+$
there is $M$ such that $s(r)\leq f(r)+M$. Equivalently,
there is $r_0 > 0$ such that $s(r)\leq f(r)$ for $r\ge r_0$.
\end{Def}

\begin{Prop}\label{ConstructionOfSublinear}
If $t_n,a_n\in\R_+$ and $t_n\to\infty$, $a_n\to 0$,
then there is an asymptotically sublinear function $s\colon \R_+\to\R_+$
such that $s(t_k)=a_k\cdot t_k$ for infinitely many $k$.
\end{Prop}
\dokaz If for each $n$ there is $k > n$ such that $a_k\cdot t_k\leq a_n\cdot t_n$, then such $s$ is easy to construct as a piecewise linear function
that is decreasing starting from some $t$.
Without loss of generality assume $a_k\cdot t_k > a_n\cdot t_n$ if $k > n$.
The idea is to pick a sequence $n(i)$ such that
the slopes $\frac{a_{n(i+1)} t_{n(i+1)}-a_{n(i)} t_{n(i)}}{t_{n(i+1)}- t_{n(i)}}$
form a decreasing sequence tending to $0$. Once it is done (easy exercise) create the graph of $s$ in a piecewise linear fashion.
\edokaz

\begin{Def}\label{HigsonSublinearMapsDef}
A continuous function $f\colon (X,d_X)\to (Y,d_Y)$ is called {\it Higson
sublinear} if for each asymptotically sublinear function $s\colon \R_+\to\R_+$
the conditions $x_n,y_n\to\infty$ and $d_X(x_n,y_n)\leq s(|x_n|)$ for all $n\ge 1$ imply
$d_Y(f(x_n),f(y_n))\to 0$.
\end{Def}

Using asymptotically sublinear functions one can introduce
the {\it sublinear coarse structure} $\EEE_L(X)$ on $X$
either following \cite{DydHof} and declaring a family of subsets of $X$ to be uniformly
bounded if and only if it refines $\{B(x,s(|x|))\}_{x\in X}$
for some asymptotically sublinear $s$.
Alternatively, one can use \cite{Roe lectures} 
and declare $E\in\EEE_L(X)$ if and only if
$E\subset \bigcup\limits_{x\in X}B(x,s(|x|))\times B(x,s(|x|))$ for some
asymptotically sublinear $s$.
The corresponding {\it Higson sublinear compactification} $h_L(X)$
has the property that all bounded maps $f\colon h_L(X)\to \R$
induce Higson sublinear maps $f\vert_X$ and every
Higson sublinear bounded $f\colon X\to \R$ extends over $h_L(X)$.
If $X$ is proper, then $\nu_L(X):=h_L(X)\setminus X$
is called the {\it sublinear Higson corona of $X$}.

\begin{Prop}\label{LipschitzImpliesHigsonSublinearF}
Suppose $f\colon X\to S^m$ is a continuous function.
If $f'\colon X\to \R^{m+1}$ defined by $f'(x)=|x|\cdot f(x)$
is asymptotically Lipschitz, then $f$ is Higson sublinear.
\end{Prop}
\dokaz Suppose $s\colon\R_+\to\R_+$ is asymptotically sublinear
and $x_n,y_n\to\infty$ so that $d(x_n,y_n)\leq s(|x_n|)$ for all $n$.
Let $\lambda, M > 0$ satisfy $|f'(x)-f'(y)|\leq \lambda\cdot d(x,y)+M$ for all $x,y\in X$
and suppose $\epsilon > 0$.
Pick $r_0 > 0$ such that $s(r)\leq \frac{\epsilon}{\lambda+1}\cdot r$
for $r > r_0$. If $n$ is sufficiently large so that $|x_n| > r$,
then
$$|f(x_n)-f(y_n)|\cdot |x_n|=||x_n|\cdot f(x_n)-|x_n|\cdot f(y_n)|=
|f'(x_n)-f'(y_n)-(|x_n|-|y_n|)\cdot f(y_n)|\leq $$
$$|f'(x_n)-f'(y_n)|+||x_n|-|y_n||\leq \lambda\cdot d(x_n,y_n)+M+d(x_n,y_n)\leq (\lambda+1)\cdot s(|x_n|)+M$$ $$\leq \epsilon \cdot |x_n|+M.$$
Therefore $|f(x_n)-f(y_n)|\leq 2\epsilon$ for large $n$.
\hfill $\blacksquare$

\begin{Rem}\label{RemOnDiffBetwHigsonAndStrongLip}
Consider $f\colon \NN\to\ S^1$ defined by $f(2n)$ being at distance $\frac{1}{\sqrt{n}}$
from $1$
and $f(2n+1)=1$ for all $n\in\NN$. Notice that $f'(n)=n\cdot f(n)$
is not asymptotically Lipschitz but $f$ is Higson sublinear.
\end{Rem}

Our next result generalizes part of Lemma 2.3 of \cite{DRS} and our proof seems simpler.
\begin{Prop}\label{CoversInducedFromSublComp}
Suppose $(X,d)$ is a discrete metric space.
If $\UUU=\{U_1,\ldots,U_k\}$ is an open cover
of $h_L(X)$, then either $X\subset U_i$ for some $i$ or there is $\epsilon > 0$ such that
$\sum\limits_{i=1}^kd(x,X\setminus U_i)\ge \epsilon\cdot |x|$
for all $x\in X$.
\end{Prop}
\dokaz By contradiction. Using \ref{ConstructionOfSublinear} find a sequence of points
$x_n\in X$ tending to infinity and an asymptotically sublinear function $s$
such that $\sum\limits_{i=1}^kd(x_n,X\setminus U_i) < s( |x_n|)$
for all $n\ge 1$.
Given that data, the sequence $\{x_n\}$ has a cluster point $x_0\in \nu_L(X)$.
As $x_0\in U_i$ for some $i$, we may as well assume $x_n\in U_1$
for all $n\ge 0$. Since $d(x_n,X\setminus U_1) < s(|x_n|)$
for all $n\ge 1$, either $X\subset U_1$ or
for each $n$ there is $y_n\in X\setminus U_1$ such that
$d(x_n,y_n) < s(|x_n|)$. Pick a continuous function
$u\colon h_L(X)\to [0,1]$ such that $u(x_0)=1$ and $u(h_L(X)\setminus U_1)=\{0\}$.
$u\vert_X$ is Higson sublinear, so $1=u(x_n)-u(y_n)\to 0$,
a contradiction.
\edokaz

\begin{Prop}\label{CharOfHigsonSublinearF}
Suppose $(X,d)$ is a discrete metric space.
A function $f\colon X\to S^m$ is Higson sublinear if and only
it is the uniform limit of functions $g\colon X\to S^m$
such that the functions $g'\colon X\to \R^{m+1}$, $g'(x)=|x|\cdot g(x)$,
are Lipschitz.
\end{Prop}
\dokaz Suppose $f\colon X\to S^m$ is Higson sublinear.
Replace $S^m$ by $\partial\Delta^{m+1}$ using a bi-Lipschitz homeomorphism. Given $\epsilon > 0$
approximate $f$ by a partition of unity on a cover of $X$ that is a restriction
of an open cover of $h_L(X)$. Such partitions of unity
induce $g\colon X\to \partial\Delta^{m+1}$ for which
$g'\colon X\to \OC(\partial\Delta^{m+1})$ is Lipschitz (see \ref{CoversInducedFromSublComp} and \ref{CharOfNiceCoversviaLip}).
Switch back to $S^m$ using \ref{ConeEquivalence}.
\hfill $\blacksquare$

\begin{Prop}\label{HigsonCoronaOfANet}
Suppose $(X,d)$ is a proper metric space.
If $X'$ is a discrete $\epsilon$-net in $X$ for some $\epsilon > 0$,
then the inclusion $i\colon X'\to X$ induces a homeomorphism
of sublinear Higson coronas $\nu_L(X')\to \nu_L(X)$.
\end{Prop}
\dokaz Notice the image of $h_L(X')$ in $h_L(X)$ must contain
$\nu_L(X)$. Otherwise there is a continuous $\alpha\colon h_L(X)\to [0,1]$
sending $i_\ast(h_L(X'))$ to $0$ and some point $x_1\in\nu_L(X)$ to $1$.
$\alpha\vert_X$ is Higson sublinear, so there is $r > 0$ such that
$\alpha(X\setminus B(x_0,r))\subset [0,1/2]$ (otherwise there are points $x_n,y_n\in X$
tending to infinity with $x_n\in X'$, $d(x_n,y_n) < \epsilon$, and $\alpha(y_n)> 1/2$
contradicting $\alpha\vert_X$ being Higson sublinear).
Since $X\setminus B(x_0,r)$ contains $x_1$ in its closure in $h_L(X)$,
$\alpha(x_1)\leq 1/2$, a contradiction. 
Thus $\nu_L(X')\to \nu_L(X)$ is surjective.
\par Suppose $x_1\ne x_2\in \nu_L(X')$.
Pick two closed neighborhoods $A_1$ of $x_1$ and $A_2$ of $x_2$ in $h_L(X')$
that are disjoint. Notice $X\cap i(A_j)=X'\cap A_j$, $j=1,2$,
so call that set $C_j$.
By \ref{CoversInducedFromSublComp} applied
to the cover $\{h_L(X')\setminus A_1,h_L(X')\setminus A_2\}$
there is $\delta > 0$ such that $d(x,C_1)+d(x,C_2)\ge \delta\cdot |x|$
for all $x\in X'$. Let us show that implies existence of $\gamma > 0$ such that
$S(x):=d(x,C_1)+d(x,C_2)\ge \gamma\cdot |x|$
for all $x\in X$.
\par Given $x\in X$ pick $x_1\in X_1$ so that $d(x,x_1) < \epsilon$.
Since $S$ is $2$-Lipschitz, $S(x)\ge S(x_1)-2d(x,x_1)\ge \delta|x_1|-2\epsilon$
and $S(x)\ge (\delta/2)|x|$ if $|x|\ge 2\epsilon(2+\delta)/\delta$.
For $x$ such that $0 < |x|\leq 2\epsilon(2+\delta)/\delta$ pick $a_1\in C_1$
and $a_2\in C_2$ so that $d(x,a_i)\leq \lambda |x|$ for some small $\lambda > 0$.
Hence $d(a_1,a_2)\leq 2\lambda|x|$ is small but the minimal distance between $C_1$ and $C_2$
is at least the discretness of $X_1$.
By  \ref{CharOfNiceCoversviaLip}
and \ref{LipschitzImpliesHigsonSublinearF}
the function $x\to \frac{d(x,C_1)}{d(x,C_1)+d(x,C_2)}$ is Higson sublinear on $X$
and extends uniquely over $h_L(X)$.
Its restriction to $X'$ vanishes on $C_1$ and is $1$ on $C_2$
implying $i(x_1)\ne i(x_2)$.
\hfill $\blacksquare$

\section{Dimension of sublinear Higson corona}

\begin{Thm}\label{CharOfDimOfHigsonSublinearCoronaForDiscrete}
Suppose $(X,d)$ is a discrete metric space.
The dimension of the sublinear Higson compactification $h_L(X)$ of $X$
is the smallest integer $m\ge 0$ with the following property:
Any norm-preserving Lipschitz function $f'\colon A\to \R^{m+1}$, $A\subset X$, extends to a 
norm-preserving Lipschitz
function $g'\colon X\to \R^{m+1}$.
\end{Thm}
\dokaz Suppose $\dim(h_L(X))\leq m$ and
$f'\colon A\subset X\to \R^{m+1}$ is a Lipschitz function induced by $f\colon A\to S^m$.
Replace $S^m$ by $\partial\Delta^{m+1}$ using  \ref{ConeEquivalence}.
Using \ref{UpperHemisphereCase} extend $f$ to $f\colon X\to \Delta^{m+1}$ inducing Lipschitz $f'$.
Find a Lipschitz $r\colon \Delta^{m+1}\to \Delta^{m+1}$
that induces a retraction $N\to \partial\Delta^{m+1}$
from a closed neighborhood $N$ of $\partial\Delta^{m+1}$ in $\Delta^{m+1}$.
Let $g_1\colon A_1=f^{-1}(N)\to \partial\Delta^{m+1}$
be defined by $g_1(x)=r(f(x))$.
 $g_1$ is Higson sublinear by \ref{LipschitzImpliesHigsonSublinearF}.
As $\dim(h_L(X))\leq m$, it extends over $h_L(X)$ to
a continuous function $g_2$. Using \ref{CharOfHigsonSublinearF} choose $g_3\colon X\to \partial\Delta^{m+1}$ with the property that the segment joining $g_2(x)$ and $g_3(x)$
in $\Delta^{m+1}$ is always contained in $N$ for all $x\in X$ and $g_3$ induces Lipschitz $g_3'\colon X\to \OC(\partial\Delta^{m+1})$.

 Put $U_1=f^{-1}(int(N))$ and $U_2=f^{-1}(\Delta^{m+1}\setminus \partial\Delta^{m+1})$.
By \ref{CoversInducedFromSublComp} and \ref{CharOfNiceCoversviaLip}, the canonical partition of unity $\{\alpha,\beta\}$ of
the cover $\{U_1\cap X,U_2\cap X\}$ gives a map
from $X$ to $\Delta^1$ inducing Lipschitz $X\to \OC(\Delta^1)$. 
Since $\alpha\cdot (r\circ f)+\beta\cdot g_3$ induces a Lipschitz function to the open cone (see \ref{ConvexityOfSpecialHigson}) and its values are in $N$, $g=r\circ (\alpha\cdot (r\circ f)+\beta\cdot g_3)$ is the required extension of $f$.

\par Suppose $\asdim_L(X)\leq m$.
Given a map $f\colon h_L(X)\to \Delta^{m+1}$
it induces an open cover $\UUU=\{U_1,\ldots,U_{m+2}\}$
of $h_L(X)$, where $U_i=f^{-1}(St(e_i))$.
By \ref{CoversInducedFromSublComp}, \ref{CharOfNiceCoversviaLip},
and \ref{ConeEquivalence} there is $g\colon h_L(X)\to\partial\Delta^{m+1}$
with the property that $g(x)\notin St(e_i)$ if $f(x)\notin St(e_i)$.
That is sufficient to conclude that $S^m$ is an absolute extensor of $h_L(X)$, i.e. $\dim(h_L(X))\leq m$. Indeed,
if $f$ is an extension of $f_A\colon A\subset h_L(X)\to \partial\Delta^{m+1}$, then $g\vert_A$ and $f$ are homotopic (via straight segments) resulting in $f_A$ being extendible over $h_L(X)$.
\hfill $\blacksquare$

\begin{Cor}\label{CharOfDimOfHigsonSublinearCoronaForProper}
Suppose $(X,d)$ is a proper metric space.
The dimension of the sublinear Higson corona $\nu_L(X)$ of $X$
is the smallest integer $m\ge 0$ with the following property:
Any norm-preserving asymptotically Lipschitz function $f'\colon A\to \R^{m+1}$, $A\subset X$, extends to a 
norm-preserving asymptotically Lipschitz
function $g'\colon X\to \R^{m+1}$.
\end{Cor}
\dokaz Pick a discrete $1$-net $X_1$ in $X$ and notice
$h_L(X_1)$ is the closure of $X_1$ in $h_L(X)$ and $h_L(X_1)=X_1\cup \nu_L(X)$ by
\ref{HigsonCoronaOfANet}. Use \ref{EquivOfLipAndAsympLipCase}. 
\hfill $\blacksquare$

\begin{Cor}\label{DimOfHigsonSublinearCoronaAndAN}
Suppose $(X,d)$ is a cocompact proper metric space that is $M$-scale connected for some $M$.
If $\asdim_{AN}(X)$ is finite, then $\asdim_{AN}(X)=\dim(\nu_L(X))$.
\end{Cor}
\dokaz As in \ref{CharOfDimOfHigsonSublinearCoronaForProper}
consider a discrete $1$-net $X_1$ in $X$ while using \ref{CharOfAsympANDimViaNormPres} to get inequality $\asdim_{AN}(X)\leq\dim(\nu_L(X))$.
The other inequality follows from \ref{IfSmIsLipschitzExtensorThenDimIsAtMostm} using the fact $S^m$
is a Lipschitz extensor if $\dim_{AN}(X)\leq m$ (see \cite{LangSch} or \cite{BDHM}).
\hfill $\blacksquare$

\section{Final comments}
 In \cite{Dr1} Dranishnikov introduced the category $\tilde \AAA$ whose objects
 are pointed proper metric spaces $X$ and morphisms are asymptotically Lipschitz
 functions $f\colon X\to Y$ such that there are constants $b,c > 0$ satisfying
 $|f(x)|\ge c\cdot |x|-b$ for all $x\in X$.
 Let us characterize $\asdim_L(X)$ in terms of $\tilde \AAA$. Namely, $\asdim_L(X)\leq n$
 if and only if $\R^{n+1}$ is an absolute extensor of $X$ in the category $\tilde\AAA$.
 
 \begin{Prop}\label{NormPreservingVsATilde}
Suppose $(X,d)$ is a pointed proper metric space and $s\colon X\to \R_+$
is a morphism of $\tilde \AAA$. If $f\colon X\to S^n$, then the following conditions are equivalent:
\begin{itemize}
\item[a.] $f'\colon X\to \R^{n+1}$, $f'(x)=|x|\cdot f(x)$, is asymptotically Lipschitz.
\item[b.] $F\colon X\to \R^{n+1}$ is asymptotically Lipschitz, where $F(x)=s(x)\cdot f(x)$.
\end{itemize}
\end{Prop}
\dokaz Suppose $\lambda\cdot |x|+M\ge s(x)\ge 2c\cdot |x|-b$ for some constants
$\lambda, M, c, b> 0$. Also, $|s(x-s(y)|\leq \lambda\cdot d(x,y)+M$.
\par 
a)$\implies$b). Assume $|f'(x)-f'(y)|\leq u\cdot d(x,y)+v$ for some $u,v > 0$.
If $|y|\leq b/c$, then $|F(x)-F(y)|\leq s(x)+s(y)\leq 
\lambda \cdot |x|+M +\lambda \cdot b/c+M\leq
\lambda \cdot d(x,y)+\lambda \cdot b/c+\lambda \cdot b/c+2M$,
so assume $|x|\ge |y| > b/c$.
Notice $f'(x)-f'(y)=|x|\cdot (f(x)-f(y))+(|x|-|y|)\cdot f(y)$,
so $|x|\cdot |f(x)-f(y)|\leq u\cdot d(x,y)+v+d(x,y)$.
Similarly, $F(x)-F(y)=s(x)\cdot (f(x)-f(y))+(s(x)-s(y))\cdot f(y)$
and $|F(x)-F(y)|\leq \frac{s(x)}{|x|}\cdot (u\cdot d(x,y)+v+d(x,y))+\lambda\cdot d(x,y)+M\leq
(\lambda+Mc/b)\cdot (u\cdot d(x,y)+d(x,y)+v)+\lambda\cdot d(x,y)+M$.

\par b)$\implies$a). By enlarging $\lambda$ and $M$, if necessary,
assume $|F(x)-F(y)|\leq \lambda d(x,y)+M$ for all $x,y\in X$.
If $|y|\leq b/c$, then $|f'(x)-f'(y)|\leq |x|+|y|\leq 
d(x,y)+2|y|\leq d(x,y)+2b/c$,
so assume $|x|\ge |y| > b/c$ and observe $s(x)\ge c\cdot |x|$.
Notice $F(x)-F(y)=s(x)\cdot (f(x)-f(y))+(s(x)-s(y))\cdot f(y)$,
so $s(x)\cdot |f(x)-f(y)|\leq 2\lambda\cdot d(x,y)+2M$.
Similarly, $f'(x)-f'(y)=|x|\cdot (f(x)-f(y))+(|x|-|y|)\cdot f(y)$,
so $|f'(x)-f'(y)|\leq \frac{|x|}{s(x)}(2\lambda\cdot d(x,y)+2M)+d(x,y)\leq
(2\lambda\cdot d(x,y)+2M)/c+d(x,y)$.
\hfill $\blacksquare$

\begin{Cor}\label{DimOfHigsonAndATilde}
Suppose $(X,d)$ is a pointed proper metric space.
The dimension of the sublinear Higson corona of $X$
is the smallest integer $n\ge 0$ such that
$\R^{n+1}$ is an absolute extensor of $X$ in the category $\tilde \AAA$.
\end{Cor}
\dokaz Suppose $\dim(\nu_L(X))\leq n$
and $F\colon A\subset X\to \R^{n+1}$ is a morphism of $\tilde\AAA$,
where $A$ is a closed subset of $X$.
Since $\R_+$ is an absolute extensor of $X$ in $\tilde\AAA$ by Theorem 3.1 of \cite{Dr1}
(p.1095),
there is a morphism $s\colon X\to \R_+$ of $\tilde\AAA$
such that $s(x)=|F(x)|$ for all $x\in A$.
Pick $f\colon A\to S^n$ such that $F(x)=s(x)\cdot f(x)$ for all $x\in A$
($f(x)$ is not uniquely defined if $s(x)=0$ in which case pick a value for $f(x)$ in
an arbitrary fashion).
Extend $f'\colon A\to \R^{n+1}$ to a norm-preserving asymptotically
Lipschitz $g'\colon X\to \R^{n+1}$ induced by $g\colon X\to S^n$
(use \ref{NormPreservingVsATilde} and \ref{CharOfDimOfHigsonSublinearCoronaForProper}). Define $G\colon X\to \R^{n+1}$ by $G(x)=s(x)\cdot g(x)$.
It is asymptotically Lipschitz by \ref{NormPreservingVsATilde}.
\par Suppose $\R^{n+1}$ is an an absolute extensor of $X$ in $\tilde\AAA$.
Pick a discrete $1$-net $X_1$ in $X$. Use \ref{NormPreservingVsATilde}
to conclude that any norm-preserving Lipschitz $A\subset X_1\to \R^{n+1}$
extends over $X_1$ to a norm-preserving Lipschitz function.
By \ref{HigsonCoronaOfANet} and \ref{CharOfDimOfHigsonSublinearCoronaForProper}, $\dim(\nu_L(X)\leq n$.
\hfill $\blacksquare$

Notice in the proof of \ref{DimOfHigsonAndATilde} we had to insist that $A$ is a closed subset of $X$. It is so because that is how absolute extensors of categories $\AAA$ and $\tilde\AAA$ are
defined in \cite{Dr1}. And that makes sense since, for $f\colon A\subset X\to Y$
to be a morphism in any of those categories, $A$ needs to be proper metric, hence closed in $X$.
However, once one leaves the language of categories, the requirement of $A$ to be closed can be dropped as shown in the next proposition.

\begin{Prop}\label{AbsExtensorsInATilde}
Suppose $(X,d_X)$ and $(Y,d_Y)$ are pointed metric spaces
and $f\colon A\subset X\to Y$ is an asymptotically Lipschitz function
such that $|f(x)|\ge c\cdot |x|-b$ for some constants $c,b > 0$.
If $R > 0$, then there is an asymptotically Lipschitz extension $g\colon B(A,R)\to Y$
of $f$ such that $|g(x)|\ge c\cdot |x|-v$ for some constant $v > 0$.
\end{Prop}
\dokaz Suppose $f$ is $(\lambda,M)$-asymptotically Lipschitz.
Pick a retraction $r\colon B(A,R)\to A$ such that $d_X(x,r(x))< R$ for all $x\in B(A,R)$
and define $g(x)=f(r(x))$ for $x\in B(A,R)$.
Using $|r(x)|\ge |x|-R$ one easily gets $|g(x)|=|f(r(x)|\ge c|r(x)|-b\ge c|x|-cR-b$.
Also, $d_Y(g(x),g(y))=d_Y(f(r(x)),f(r(y)))\leq \lambda\cdot d_X(r(x),r(y))+M\leq
\lambda\cdot d_X(x,y)+2\lambda R+M$.
\hfill $\blacksquare$

Using \ref{AbsExtensorsInATilde} one can extend any asymptotically Lipschitz function
from $A$ to its closure since the ball $B(A,R)$ contains that closure.

\end{document}